\title{ ~~~~~\\ Asymptotically exact heuristics for (near) primitive
roots. II  }
\author{Pieter Moree}
\documentclass[12pt]{article}
\usepackage{amssymb, latexsym, amsfonts}
\def\@ptsize{2}
\newtheorem{Thm}{Theorem}
\newtheorem{Lem}{Lemma}

\newtheorem{Cor}{Corollary}

\newcommand{\qed}{\hfill $\Box$}

\begin{document}
\date{}
\maketitle
{\def\thefootnote{}
\footnote{{\it Mathematics Subject Classification (2001)}.
Primary 11R45; Secondary 11A07}}
\begin{abstract}
\noindent Let $g\in\mathbb Q^*$ be a rational number.
Let $N_{g,t}(x)$ denote the number of primes $p\le x$
for which the subgroup of~$(\mathbb Z/p\mathbb Z)^*$ generated by
$g{\rm ~mod~}p$ is of residual index $t$. In \cite{MO} an heuristic
for $N_{g,t}(x)$ was set up, under the
assumption of the Generalized
Riemann Hypothesis (GRH), and shown to be asymptotically
exact.
In this paper we provide an alternative and rather shorter
proof of this result.\\
\indent Keywords:~heuristic, residual index, natural density,
primitive root.
\end{abstract}
\section{Introduction}
Let $g\in \mathbb Q\backslash \{-1,0,1\}$ and
$t\ge 1$ be an
arbitrary natural number. We write $g=\pm g_0^h$, where $g_0>0$ is
not an exact power of a rational and $h\in \mathbb Z_{\ge 1}$.
Every prime $p$ in this paper is (mostly
tacitly) assumed to be odd and satisfy ord$_p(g)=0$, e.g. $\pi(x;t,1)$ denotes
the number of odd primes $p\le x$ with $p\equiv 1({\rm mod~}t)$ and
ord$_p(g)=0$. We define
$r_g(p)=[(\mathbb Z/p\mathbb Z)^*:\langle g ~{\rm mod~}p\rangle]$
and say
that
$r_g(p)$ is the {\it residual index mod} $p$ of $g$. For an arbitrary natural
number $t$ we consider the set
$N_{g,t}$ of primes $p$ satisfying $r_g(p)=t$ and
let $N_{g,t}(x)$ denote the number of primes $p\le x$ in $N_{g,t}$.
Notice that $N_{g,1}$ is the set of primes $p$ such that $g$ is a primitive
root mod $p$.
In this paper we are interested in the behaviour of
$N_{g,t}(x)$ as $x$ tends to infinity. Our heuristic
approach to $N_{g,t}(x)$
will be entirely based on a heuristic approach to $R_{g,t}(x)$, which is
defined as the number of primes $p\le x$ with $t|r_g(p)$.\\
\indent Let $p$ be a prime with $p\equiv 1({\rm mod~}t)$. Note that the
density of elements $\gamma
\in (\mathbb Z/p\mathbb Z)^*$ such that 
$r_{\gamma}(p)=t$ is $\varphi((p-1)/t)/(p-1)$,
where $\varphi$ denotes Euler's totient function. Thus naively one might
expect that 
\begin{equation}
\label{simpelzeg}
N_{g,t}(x)\sim \sum_{p\le x\atop p\equiv 1({\rm mod~}t)}
{\varphi((p-1)/t)\over p-1}~~(x\rightarrow \infty).
\end{equation}
Despite the fact that for arbitrary $C>1$ and $m>
{\rm exp}(4\sqrt{\log x\log \log x})$, we have by \cite[Theorem 7]{MO},
$${1\over m}\sum_{g\le m}N_{g,t}(x)=
\sum_{p\le x\atop p\equiv 1({\rm mod~}t)}{\varphi((p-1)/t)\over p-1}
+O\left({x\over \log^C x}\right),$$
and thus the naive heuristic holds true on average, it can be shown 
on GRH, that (\ref{simpelzeg}) is false for many $g$.
In \cite{MO}, however, the following modified heuristic involving
a function $w_{g,t}(p)\in \{0,1,2\}$ was introduced and
shown to be asymptotically exact under GRH. (We stipulate that
zero multiplied by something not well-defined equals zero. By $(a,b)$ we denote
the greatest common divisor of $a$ and $b$.)
\begin{Thm}
\label{main} {\rm \cite{MO}}.
Let $g\in \mathbb Q\backslash \{-1,0,1\}$
 and $t\ge 1$ be an arbitrary integer.
Write $g=\pm g_0^h$, where $g_0\in \mathbb Q$ is positive
and not an exact power of a rational and $h\in \mathbb Z_{\ge 1}$.
Let $d(g_0)$ denote the discriminant of $\mathbb Q(\sqrt{g_0})$.
Let $2^e||h$ and $2^{\tau}
\mid \mid t$. Put $h_t=h/(h,t)$ and
$$\epsilon_1=\cases{0 &if $\tau<e$; \cr -1 &if $\tau=e$;\cr
1 &if $\tau>e$.}$$ If $g>0$, $p\equiv 1({\rm mod~}t)$
and $((p-1)/t,h_t)=1$, then put
$$w_{g,t}(p)= 1+{\epsilon_1\over 2}\{1+(-1)^{p-1\over 2^e}\}\left({d(g_0)\over
p}\right).$$
If $g<0$, $2\nmid h_t$, $p\equiv 1({\rm mod~lcm}(2^{e+1},t))$ and
$$((p-1)/t,h_t)=1,$$ then put
$$w_{g,t}(p)=
1+\epsilon_1(-1)^{p-1\over 2^{e+1}}\left({d(g_0)\over p}\right).$$
If $g<0$, $2|h_t$, $p\equiv 1({\rm mod~}2t)$ and $((p-1)/2t,h_t)=1$,
then
put $w_{g,t}(p)=2$. In all other cases put $w_{g,t}(p)=0$.\\
\indent Under GRH we have
\begin{equation}
\label{whatitisabout}
N_{g,t}(x)=(h,t)\sum_{p\le x\atop p\equiv 1({\rm mod~}t)}
w_{g,t}(p){\varphi((p-1)/t)\over p-1}+O\left({x\log \log x\over
\log^{2}x}\right).
\end{equation}
\end{Thm}
The purpose of this note is to give a much shorter proof of
Theorem 1 than the one given in \cite{MO}. The asymptotically
exact heuristics developed here for $R_{g,t}(x)$ have further
applications, for example in the study of exact heuristics
for divisors of recurrences of second order \cite{MOf}.

\vfil\eject

\section{Results of Hooley and Wagstaff}
In this section we briefly recall the approach of Hooley and Wagstaff in
estimating $N_{g,t}(x)$, cf. \cite{WA}; it is analogous to
Hooley's \cite{HO} break through attack on $N_{g,1}(x)$, the primitive
root counting function.
The basic observation is that $t|r_g(p)$ iff
$p$ splits completely in the splitting field
$F_t=\mathbb Q(\zeta_t, g^{1/t})$ of the polynomial $X^t-g$ over
$\mathbb Q$, where $\zeta_t=\exp(2\pi i/t)$.
As a consequence (Corollary \ref{pitconsequence}) of the prime ideal
 theorem, a special case of both the Frobenius and the Chebotarev
 density theorem,
 the set of these primes
has natural density $1/[F_t:\mathbb Q]$.
The primes in $N_{g,t}$ are
those that do not split completely in any of the
fields $F_{kt}$ with $k>1$.
A standard inclusion-exclusion argument readily yields the heuristic
value
\begin{equation}
\label{graden}
\sum_{k=1}^\infty {\mu(k)\over [F_{kt}:\mathbb Q]},
\end{equation}
for the natural density of  the set $N_{g,t}$ (provided that $N_{g,t}$
has indeed a natural  density).
The sum (\ref{graden}) converges whenever $g$ is
different from $\pm 1$,
since in that case $[F_{kt}:\mathbb Q]$ differs from its `approximate value'
$kt\varphi(kt)$
by a factor that is easily bounded in terms of $g$,
cf. Lemma \ref{wagstaffdegree}.
In fact, we obtain an upper bound for the upper density of
the set $N_{g,t}$
in this way.
In order to turn this heuristic argument for
$N_{g,1}$ into a proof, Hooley
employed estimates for the remainder term in the prime number theorem
for the fields $F_{k}$
that are currently only known to hold under GRH. Hooley's arguments
are easily extended to $N_{g,t}$ and result in the following theorem.
\begin{Thm} {\rm (GRH)}.
\label{sammie}
Let $g\in \mathbb Q\backslash\{-1,0,1\}$ and
${\rm Li}(x)=\int_2^x{\rm dt}/\log t$.
Then
$$N_{g,t}(x)=A(g,t){\rm Li}(x)+O\left({x\log \log x\over
\log^2 x}\right),$$
where
$$A(g,t):=\sum_{k=1}^{\infty}{\mu(k)\over [\mathbb Q(\zeta_{kt},g^{1/kt}):
\mathbb Q]}.$$
\end{Thm}
Now $A(g,t)$ can be expressed as a linear combination of sums of the
form
\begin{equation}
\label{wagstaffsum}
S(h,t,m):=\sum_{k=1\atop m|kt}^{\infty}
{\mu(k)(kt,h)\over kt\varphi(kt)}.
\end{equation}
Each sum $S(h,t,m)$ can be written as an Euler product and in fact
is a rational multiple of $A$, the Artin constant (cf. Lemma 2.1
and Theorem 2.2 of \cite{WA}).

\vfil\eject
\section{Proof preview}
In this section we sketch the proof of Theorem \ref{main}
given in \cite{MO} and give a preview of the proof of
Theorem \ref{main} to be given in section \ref{newproof},
making also clear the advantages of the new proof over
the old one.\\
\indent We begin by sketching the proof of Theorem \ref{main}
given in \cite{MO}. Let $G$ be a finite cyclic group of
order $n$, $t|n$ and $\gamma \in G$.
We put $f'_{\gamma,t}(G)=1$ if $[G:\langle \gamma\rangle]=t$
and $f'_{\gamma,t}(G)=0$ otherwise. It is easy to show that
\begin{equation}
\label{ouderekaridentiteit}
f'_{\gamma,t}(G)={\varphi(n/t)\over n}\sum_{d|n}
{\mu(d/(d,t))\over \varphi(d/(d,t))}
\sum_{{\rm ord}\chi=d}\chi(\gamma),
\end{equation}
where the sum is over all
multiplicative characters of $G$ of order $d$. Let us
for simplicity assume that $g>0$. We notice that if
$r_g(p)=t$, then $((p-1)/t,h_t)=1$. Thus
\begin{equation}
\label{ouwetje}
N_{g,t}(x)=\sum_{p\le x\atop ({p-1\over t},h_t)=1}f'_{g,t}
((\mathbb Z/p\mathbb Z)^*).
\end{equation}
In the latter sum we calculate
the contribution of the characters $\chi(g_0^h)$ with
$\chi$ of order $d$ with
$d|h$ (the `linear' contribution) and those $\chi$ of
order dividing $2h$ but not $h$ (the `quadratic' contribution).
For those in the linear contribution we have $\chi(g)=1$
and for those in the quadratic we have $\chi(g)=(d(g_0)/p)$.
The linear contribution turns out to be equated with the naive
heuristic approach and the linear together with the quadratic
contribution with a more subtle heuristic approach based
on having a
priori knowledge of $(d(g_0)/p)$.
Working out the contributions of the relevant
characters one obtains the sum appearing in (\ref{whatitisabout}).
In doing so crucial use of the condition
$((p-1)/t,h_t)=1$ is made (note that (\ref{ouwetje})
is also
valid when the condition $((p-1)/t,h_t)=1$ is dropped).
The sum in (\ref{whatitisabout}) can be unconditionally evaluated
with error term $O(x\log ^{-C}x)$, with $C>1$ arbitrary.
It turns out to be most convenient to do so in terms of
the sums $S(h,t,m)$ defined
in (\ref{wagstaffsum}). This allows them to be compared,
under GRH,
with Wagstaff's evaluation of $N_{g,t}(x)$
(\cite[Theorem 2.2]{WA}). The latter
evaluation requires several arithmetically complicated cases
to be distinguished. This makes the comparison process
rather involved. Theorem 1 then
follows on noticing that in all cases we
have equality up to the required error term.\\
\indent In the proof of Theorem 1 given here we
start out by considering
$f_{\gamma,t}(G)$, which is defined as $f_{\gamma,t}(G)=1$ if
$t|[G:\langle \gamma\rangle]$ and $f_{\gamma,t}(G)=0$ otherwise. The
analog of (\ref{ouderekaridentiteit}) is (\ref{startid}),
which, being arithmetically easier, is less complicated to work with.
Proceeding as before we end up concluding that the
`linear' and `quadratic' contribution taken together yield an
asymptotically exact heuristic for $R_{g,t}(x)$. The comparison
process is easier here and only requires evaluating
$[\mathbb Q(\zeta_{kt},g^{1/kt}):\mathbb Q]$ instead of $A(g,t)$.
By means of (\ref{iniexi}) the heuristic for $R_{g,t}(x)$ is
easily `pushed through' to a heuristic for $N_{g,t}(x)$. In this
approach one bypasses the use of the sums $S(h,t,m)$ and the
condition $((p-1)/t,h_t)=1$ comes up naturally (Lemma
\ref{inclucond} with $G=(\mathbb Z/p\mathbb Z)^*$ and 
thus $n=p-1$), instead
of in a somewhat ad hoc way.

\vfil\eject
\section{On $t$-divisible residual indices}
\label{tdeelbaar}
Let $g\in \mathbb Q \backslash \{-1,0,1\}$. Let
$R_{g,t}$ be the set of primes $p$ with $t|r_g(p)$ and
$R_{g,t}(x)$ the number of primes $p\le x$ in $R_{g,t}$.\\
\indent We study the set $R_{g,t}$ by two different methods. On the
one hand by characters of $(\mathbb Z/p\mathbb Z)^*$, where $p$ runs
over the primes and on the other hand by algebraic
and analytic number theory.
 Thus
in Lemma \ref{divchar} we set up a character identity.
In Lemma \ref{linquad} the sum, $M_{g,t}(x)$, of
the contributions of the `linear' and `quadratic'
characters to $R_{g,t}(x)$ is evaluated, invoking Lemma
\ref{charactersum}. The asymptotic behaviour of $M_{g,t}(x)$ is
determined
(in Lemma \ref{mtxasymptotic}) and compared with the asymptotic
behaviour
 of $R_{g,t}(x)$
as resulting from  number theory (given
in Lemma \ref{algnotheory}). This comparison then shows that the total contribution
of the characters of third and higher order is of lower order than
the contribution of the characters of first and second order
(Theorem
\ref{surprise1}).

\begin{Lem}
\label{divchar}
Let $G$ be a finite cyclic group of order $n$, $t|n$ and $\gamma\in G$.
Put $f_{\gamma,t}(G)=1$ if $t|[G:\langle \gamma\rangle]$ and
$f_{\gamma,t}(G)=0$ otherwise.
Then
\begin{equation}
\label{startid}
f_{\gamma,t}(G)={1\over t}\sum_{d|t}\sum_{{\rm ord}\chi=d}\chi(\gamma),
\end{equation}
where the inner summation is over the multiplicative characters on $G$
having order precisely $d$.
\end{Lem}
{\it Proof}. First consider the case where $t$ is squarefree. On noting
that $\sum_{{\rm ord}\chi=d}\chi(\gamma)$ is multiplicative in $d$, we
obtain
$$f_{\gamma,t}(G)={1\over t}\prod_{p|t}\left(1+\sum_{{\rm ord}\chi=p}\chi(\gamma)
\right).$$
If $p|[G:\langle \gamma\rangle]$, then
$\sum_{{\rm ord}\chi=p}\chi(\gamma)=p-1$
and hence $f_{\gamma,t}(G)=1$ if $t|[G:\langle \gamma\rangle]$.
On the other hand,
if $t\nmid [G:\langle \gamma\rangle]$, then there is a prime $q$ such that
$q|t$ and $q\nmid [G:\langle \gamma\rangle]$. Then
$\sum_{{\rm ord}\chi=q}\chi(\gamma)=-1$ and $f_{\gamma,t}(G)=0$.
 The general case
is not so immediate, but easily dealt with on using Proposition 5
of \cite{MO}. \qed\\

\noindent Lemma \ref{divchar} and its proof can also be formulated in
terms of Ramanujan sums
$$c_d(n):=\sum_{1\le k\le d\atop (k,d)=1}e^{2\pi ikn\over d}.$$
\begin{Lem}
\label{ramasum}
We have $f_{\gamma,t}(G)={1\over t}\sum_{d|t}c_d([G:\langle
\gamma \rangle])$.
\end{Lem}
{\it Proof}. Suppose $e|r$, then
by \cite[p. 79]{Mac} the sum
$\sum_{d|r}c_d(e)$ equals $r$ if $e=r$ and zero otherwise.
Now let $e$ and $r$ be arbitrary. On noticing that
$\sum_{d|r}c_d(e)=\sum_{d|r}c_d((e,r))$ and invoking the latter
result, we see that
$\sum_{d|r}c_d(e)$ equals $r$ if $r|e$ and zero otherwise.
On putting $r=t$ and $e=[G:\langle \gamma\rangle]$ the result follows. \qed\\

\noindent {\tt Remark}. By M\"obius inversion it follows from Lemma \ref{divchar}
and Lemma \ref{ramasum} that
$$\sum_{{\rm ord}\chi=d}\chi(\gamma)
=\sum_{\delta |d}\delta f_{\gamma,\delta }(G)\mu({d\over \delta})
=c_d([G:\langle \gamma\rangle ]).$$
\indent For $p$ prime put $G=(\mathbb Z/p\mathbb Z)^*$
and $\psi_d(p)=\sum_{{\rm ord}\chi=d}\chi(g)$.
\begin{Lem}
\label{charactersum}
Adopt the notations and assumptions of Theorem 1.
Put
$$\epsilon_2=\cases{0 &if $\tau\le
e$;\cr 1&if $\tau> e$. }$$
Assume that $p\equiv 1({\rm mod~}(h,t))$.
If $g>0$, then $\sum_{d|(h,t)}\psi_{d}(p)=(h,t)$.\\
If $g<0$, then
$$\sum_{d|(h,t)}\psi_{d}(p)=\cases{(h,t) &if $p\equiv 1({\rm
mod~}2(h,t))$;
\cr 0 &otherwise.}$$
Assume that $p\equiv 1({\rm mod~}(2h,t))$. If $g>0$, then
$$\sum_{d|(2h,t)\atop d\nmid h}\psi_{d}(p)=\epsilon_2 ({d(g_0)\over p})
(h,t).$$
If $g<0$, then
$$\sum_{d|(2h,t)\atop d\nmid h}\psi_{d}(p)=\epsilon_2 (-1)^{p-1\over
2^{e+1}}
({d(g_0)\over p})(h,t).$$
\end{Lem}
{\it Proof}. Straightforward on using Lemma 15 of \cite{MO} to
evaluate $\psi_d(p)$ in each of the four cases, cf. the proof
of Lemma 16 of \cite{MO}. \qed\\

\noindent Put $$L_{g,t}(x)={1\over t}
\sum_{p\le x\atop p\equiv 1({\rm mod~}t)}\sum_{d|(h,t)}
\psi_d(p) {\rm ~and~} Q_{g,t}(x)=
{1\over t}\sum_{p\le x\atop p\equiv 1({\rm mod~}t)}
\sum_{d|(2h,t)\atop d\nmid h}\psi_{d}(p).$$
Using Lemma \ref{divchar}, the definitions of $\psi_d(p),~L_{g,t}(x)$
and $Q_{g,t}(x)$, we conclude that
$$R_{g,t}(x)=L_{g,t}(x)+Q_{g,t}(x)+{1\over t}\sum_{p\le x\atop
p\equiv 1({\rm mod~}t)}\sum_{d|t\atop d\nmid 2h}\psi_d(p).$$
Although we are interested only in $M_{g,t}(x):=L_{g,t}(x)+Q_{g,t}(x)$, $L_{g,t}(x)$ and
$Q_{g,t}(x)$ are of rather different nature and hence we are forced to
consider them by themselves. Roughly speaking $L_{g,t}(x)$ gives the
contribution of those characters such that $\chi(g)=1$ for all
characters $\chi$ having the  same order and $Q_{g,t}(x)$ of those
such that $\chi(g)$ reduces to a quadratic character for all $\chi$
having the same order.\\
\indent Let $K$ be an arbitrary algebraic number field
and let $P_K(x)$ denote the number of rational primes $p\le x$ that
split completely in $K$.
On using the previous lemma and noticing
that
$$\sum_{p\le x\atop p\equiv 1({\rm mod~}t)}({d(g_0)\over p})
=2P_{\mathbb Q(\zeta_t,\sqrt{g_0})}(x)-\pi(x;t,1),$$
one finds after some computation:
\begin{Lem}
\label{linquad} Define $M_{g,t}(x):=L_{g,t}(x)+Q_{g,t}(x)$ and $t_h=t/(t,h)$.
If $g>0$, then
$$M_{g,t}(x)=\cases{\pi(x;t,1)/t_h &if $\tau\le e$;\cr
2P_{\mathbb Q (\zeta_t,\sqrt{g_0})}(x)/t_h &if $\tau>e$.}$$
If $g<0$, then
$$M_{g,t}(x)=\cases{\pi(x;2t,1)/t_h &if $\tau\le e$;\cr
\{4P_{\mathbb Q(\zeta_{2t},\sqrt{g_0})}(x)-2P_{\mathbb Q (\zeta_t,\sqrt{g_0})}(x)
+2\sum_{{p\le x\atop p\equiv 1({\rm mod~}t)}\atop p\not\equiv
1({\rm mod~}2t)}1 \}/t_h&if $\tau=e+1$;\cr
2P_{\mathbb Q(\zeta_t,\sqrt{g_0})}(x)/t_h &if $\tau>e+1$.}$$
\end{Lem}
We will use Lemma \ref{linquad} to deduce Lemma \ref{mtxasymptotic},
which gives the asymptotic behaviour of $M_{g,t}(x)$. In order to do
so we need a result due to Siegel and Walfisz
and the prime ideal theorem (due to Landau).
\begin{Lem}
\label{SW} {\rm \cite[Satz 4.8.3]{P}}.
Let $C$ be a fixed real number. Then the estimate
$$\pi(x;d,a):=\sum_{p\le x\atop p\equiv a({\rm mod~}d)}1
={{\rm Li}(x)\over \varphi(d)}
+O(xe^{-c_1\sqrt{\log x}})$$
holds uniformly for all integers $a$ and $d$ such that
$(a,d)=1$ and $1\le d\le \log ^C x$, with $c_1$ some
positive constant.
\end{Lem}
\begin{Lem}
\label{pit} {\rm \cite{LA}}.
Let $K$ be an algebraic number field.
Let $\pi_K(x)$ denote the number of prime ideals of norm
at most $x$. There exists $c_2>0$ such that
$$\pi_K(x)={\rm Li}(x)+O(xe^{-c_2\sqrt{\log x}}).$$
\end{Lem}
\begin{Cor}
\label{pitconsequence}
Let $P_K(x)$ denote the number of rational primes $p\le x$ that
split completely in the number field $K$.
If $K$ is normal, then
$$P_K(x)={{\rm Li}(x)\over [K:\mathbb Q]}+O(xe^{-c_2\sqrt{\log x}}).$$
\end{Cor}
{\tt Remark}. A more complicated but sharper error term was obtained
by Mitsui \cite{MI}. If the Riemann hypothesis holds for the
Dedekind zeta-function $\zeta_K(s)$, then it can be shown \cite{LAN} that
the error is of
order $O(\sqrt{x}\log(d(K)x^{[K:\mathbb Q]}))$, where $d(K)$
denotes the absolute value of the discriminant of $K$.\\

\noindent Also we need an explicit evaluation of
the field degree $[\mathbb Q(\zeta_t,g^{1/t}):\mathbb Q]$, which
is given in the next lemma.
\begin{Lem}
\label{wagstaffdegree}
Write $t_h=t/(t,h)$ and
$[\mathbb Q(\zeta_t,g^{1/t}):\mathbb Q]=\varphi(t)t_h/\nu$.
If $g>0$, we have $\nu=2$ if $t_h$ is even and $d(g_0)|t$;
otherwise $\nu=1$. Now suppose $g<0$. If $t$ is odd,
then $\nu=1$. If $t$ is even and $t_h$ is odd, then
$\nu=1/2$. If $t$ is even and $t_h\equiv 2({\rm mod~}4)$, then
$$\nu=\cases{2 &if $d(g_0)\nmid t$ {\rm ~and~}$d(g_0)|2t$;\cr
1 &otherwise.}$$
If $t$ is even and $4|t_h$, then $\nu=2$ if
$d(g_0)|t$ and $\nu=1$  otherwise.
\end{Lem}
{\it Proof}. This is \cite[Proposition 4.1]{WA}, with
the condition $t\equiv 2({\rm mod~}4)$ and $d(-g_0)|t$
or $t\equiv 4({\rm mod~}8)$ and $d(2g_0)|t$ replaced
by the equivalent condition
 $d(g_0)\nmid t$ and $d(g_0)|2t$. \qed\\

\noindent With the three latter results in hand, we can now
prove the following lemma.
\begin{Lem}
\label{mtxasymptotic}
Let $C$ be a fixed real number. Then for some $c_3>0$ the estimate
$$M_{g,t}(x)={{\rm Li}(x)\over [\mathbb Q(\zeta_t,g^{1/t}):\mathbb Q]}
+O(xe^{-c_3\sqrt{\log x}})$$
holds uniformly for all $t$ with $1\le t\le (\log^C x)/d(g_0)$.
\end{Lem}
{\it Proof}. The primes $p$ that
split completely in $\mathbb Q(\zeta_t,\sqrt{g_0})$ are precisely the
primes $p$ satisfying $p\equiv 1({\rm mod~}t)$ and $(d(g_0)/p)=1$.
By the law of quadratic reciprocity these primes are precisely those
in a union of residue classes of modulus at most $4d(g_0)t$. This
means we can invoke Lemma \ref{SW}.
The natural density of the primes that split completely
in $\mathbb Q(\zeta_t,\sqrt{g_0})$ is given by Lemma
\ref{pitconsequence} as
$1/[\mathbb Q(\zeta_t,\sqrt{g_0}):\mathbb Q]$. The field degree
$[\mathbb Q(\zeta_t,\sqrt{g_0}):\mathbb Q]$ is well-known to be
$2\varphi(t)$ if $d(g_0)\nmid t$ and $\varphi(t)$ otherwise.
We obtain the assertion of the lemma with
${\rm Li}(x)/[\mathbb Q(\zeta_t,g^{1/t}):\mathbb Q]$ replaced
by ${\rm Li}(x)/c_{t,g}$, where $c_{t,g}$ is an easily explicitly
evaluated constant. On comparing the values of $c_{t,g}$ with
those of $[\mathbb Q(\zeta_t,g^{1/t}):\mathbb Q]$ as given in
Lemma \ref{wagstaffdegree}, it is seen
that $c_{t,g}=[\mathbb Q(\zeta_t,g^{1/t}):\mathbb Q]$. \qed\\

\noindent The primes in $R_{g,t}$ are easily characterized algebraically.
They are precisely the primes $p$ satisfying $p\equiv 1({\rm mod}~t)$
and $g^{(p-1)/t}\equiv 1({\rm mod~}p)$. But these are precisely the primes
splitting completely in $\mathbb Q(\zeta_t,g^{1/t})$ and thus
$R_{g,t}(x)=P_{\mathbb Q(\zeta_t,g^{1/t})}(x)$. By Corollary
\ref{pitconsequence} the following result then follows.
\begin{Lem}
\label{algnotheory}
There exists $c_4>0$ such that
$$R_{g,t}(x)={{\rm Li}(x)\over [\mathbb Q(\zeta_t,g^{1/t}):\mathbb Q]}
+O(xe^{-c_4\sqrt{\log x}}).$$
\end{Lem}
\indent Comparison of Lemma \ref{mtxasymptotic} and Lemma \ref{algnotheory}
yields:
\begin{Thm}
\label{surprise1}
We have $R_{g,t}(x)=M_{g,t}(x)+O(xe^{-c_5\sqrt{\log x}})$, for
some $c_5>0$.
\end{Thm}
This theorem can be loosely phrased as stating that
only the contributions of the `linear' and `quadratic'
characters are responsible for the
asymptotic behaviour of $R_{g,t}(x)$. That $M_{g,t}(x)$ and $R_{g,t}(x)$
are so closely related comes perhaps as a surprise, but in the next 
subsection we give a heuristic approach to $R_{g,t}(x)$  that will yield
$M_{g,t}(x)$ as outcome.

\subsection{Heuristic approach to $R_{g,t}(x)$}
Let us first consider the case $g=g_0^h$.
Then $g$ is a priori in $G^h$ (with
$G=(\mathbb Z/p\mathbb Z)^*$). We are interested in the
case where $g$
satisfies
$t|[G:\langle g \rangle]$. Note
that if $t|[G:\langle g\rangle]$, then $p\equiv 1({\rm mod~}t)$.
If $t|p-1$, then the elements of residual index $t$ are all in
$G^t$. The probability of finding $g$, given
our a priori knowledge, in $G^t$ equals
$|G^h\cap G^t|/|G^h|$. The latter quotient is the density of elements
in $G^h$ having residual index divisible by $t$ and is easily computable,
also in the case where $G$ is an arbitrary cyclic group. 
\begin{Lem}
\label{intersectie}
Let $G$ be a finite cyclic group of order $n$ and let $t$ and $h$ be
arbitrary with $t|n$. Then
$$\rho_{1,*,t}(G):={|G^h\cap G^t|\over |G^h|}={(t,h)\over t}={1\over t_h}.$$
\end{Lem}
Heuristically we might expect that
$R_{g,t}(x)$ behaves as $\sum_{p\le x\atop p\equiv 1({\rm mod~}t)}
\rho_{1,*,t}((\mathbb Z/p\mathbb Z)^*)$, that is as
$\pi(x;t,1)/t_h$. Indeed, by Lemma
\ref{linquad} and Theorem \ref{surprise1} it does, except when
$\tau>e$ and $d(g_0)|t$.
Hence let us try to refine this heuristic.
Suppose we know the value of
the Legendre symbol $(d(g_0)/p)$. This improves our a priori knowledge
and leads one to alter our group theoretical quotient. Let $\gamma$ be
a generator of $G$, thus $G=\langle \gamma \rangle$. Let $n$ be
the order of $G$. If $t|n$ we make the definitions
$$\rho_{1,1,t}(G):={|\{\gamma^{({\rm even})h}\}\cap G^t|\over
|\{\gamma^{({\rm even})h}\}|}\left( =
{|G^{2h}\cap G^t|\over |G^{2h}|}\right) {\rm ~and~}
\rho_{1,-1,t}(G):={|\{\gamma^{({\rm odd})h}\}\cap G^t|\over
|\{\gamma^{({\rm odd})h}\}|}.$$
If $t\nmid n$ we put $\rho_{1,1,t}(G)=\rho_{1,-1,t}(G)=0$. Notice that if
$(d(g_0)/p)=1$, then the reduction of $g({\rm mod~}p)$ is in
$\{\gamma^{(even)h}\}$, otherwise it is in $\{\gamma^{(odd)h}\}$.
We expect that a better heuristic for $R_{g,t}(x)$
is $H_{g,t}(x):=\sum_{p\le x}\rho_{1,(d(g_0)/p),t}((\mathbb
Z/p\mathbb Z)^*)$.
Using Lemma \ref{intersectie} one deduces
$$\rho_{1,1,t}(G)={(2h,t)\over t}{\rm ~and~}\rho_{1,-1,t}(G)=
{2(h,t)\over t}-{(2h,t)\over t}.$$
In case $\tau\le e$ this reduces to $\rho_{1,1,t}(G)=\rho_{1,-1,t}(G)=t_h^{-1}$ and
hence the naive heuristic yielding (as before),
$H_{g,t}(x)=\pi(x;t,1)/t_h$. If $\tau>e$, then
$\rho_{1,1,t}(G)=2/t_h$ and $\rho_{1,-1,t}(G)=0$, yielding
$H_{g,t}(x)=2P_{\mathbb Q(\zeta_t,\sqrt{g_0})}(x)/t_h$.
By Lemma \ref{linquad} we conclude that
$H_{g,t}(x)=M_{g,t}(x)$ in case $g>0$.\\
\indent Now suppose $g<0$ (hence $g=-g_0^h$). We assume that
$n$ is even and denote by $-1$ the unique element
of order 2 in $G$. The analog of Lemma \ref{intersectie} reads
\begin{Lem}
\label{intersectie1}
Let $G$ be a finite cyclic group of even order $n$ and let $t$ and $h$ be
arbitrary with $t|n$. Then
$$\rho_{-1,*,t}(G):={|\mbox{--}G^h\cap G^t|\over |\mbox{--} G^h|}=
\cases{0 &if ord$_2(n)=\tau$ {\rm ~and~} $\tau\le e$;\cr t_h^{-1} &otherwise.}$$
\end{Lem}
If $t\nmid n$ define 
$\rho_{-1,1,t}(G)=\rho_{-1,-1,t}(G)=0$.
If $t|n$ we make the definitions
$$\rho_{-1,1,t}(G)={|\{\mbox{--}\gamma^{({\rm even})h}\}\cap G^t|\over
|\{\mbox{--}\gamma^{({\rm even})h}\}|} {\rm ~and~}
\rho_{-1,-1,t}(G)={|\{\mbox{--}\gamma^{({\rm odd})h}\}\cap G^t|\over
|\{\mbox{--}\gamma^{({\rm odd})h}\}|}.$$
We consider how good
$H_{g,t}(x):=\sum_{p\le x}{\rho}_{-1,(d(g_0)/p),t}((\mathbb
Z/p\mathbb Z)^*)$ is as a heuristic for $R_{g,t}(x)$.
To that end we evaluate $\rho_{-1,1,t}(G)$ and $\rho_{-1,-1,t}(G)$ first.
\begin{Lem}
\label{kwadstat}
Suppose $G$ is a cyclic group of even order $n$ and $t|n$. Then
$$\rho_{-1,1,t}(G)=\cases{0 &if  ord$_2(n)=\tau$ and $\tau\le e+1$;\cr
(2h,t)/t &otherwise.}$$
Furthermore,
$$\rho_{-1,-1,t}(G)=\cases{0 &if ord$_2(n)=\tau$ and $\tau\ne e+1$;\cr
0 &if ord$_2(n)\ge \tau+1$ and $\tau\ge e+1$;\cr
(2h,t)/t &otherwise.}$$
\end{Lem}
{\it Proof}. Let us consider the more difficult case of
evaluating $\rho_{-1,-1,t}(G)$. The intersection 
$\{\mbox{--}\gamma^{({\rm odd})h}\}
\cap G^t$ consists of those elements $\gamma^{\alpha}$ with
$1\le \alpha\le n$ satisfying both
$\alpha\equiv n/2+(h,n) ({\rm mod~}(2h,n))$ and
$\alpha\equiv 0({\rm mod~}t)$. The intersection is thus
empty iff $n/2+(h,n)\not\equiv 0 ({\rm mod~}(2h,t))$.
On using that $(2h,t)$ divides both $n$ and $2(h,n)$ one infers
that $(2h,t)\nmid n/2$ and $(2h,t)\nmid (h,n)$ implies that
$(2h,t)$ divides $n/2+(h,n)$. Thus the intersection is empty
iff either $(2h,t)|n/2$ and $(2h,t)\nmid (h,n)$
or $(2h,t)\nmid n/2$ and $(2h,t)|(h,n)$. Since $(2h,t)|2(h,n)$
we have that $(2h,t)|(h,n)$ iff ord$_2((2h,t))\le {\rm ord}_2((h,n))$.
Similarly $(2h,t)|n/2$ iff ord$_2((2h,t))\le {\rm ord}_2(n)-1$.
Recalling that ord$_2(n)\ge \tau$ (by assumption),
ord$_2(h)=e$ and ord$_2(t)=\tau$, we deduce that the intersection is
empty iff either ord$_2(n)=\tau$  and $\tau\ne e+1$ or
ord$_2(n)\ge \tau+1$ and $\tau\ge e+1$.\\
\indent If the intersection is non-empty, then it
consists of $n/{\rm lcm}((2h,n),t)$,
that is $n(2h,t)(2h,n)^{-1}t^{-1}$ elements, whereas
$\{-\gamma^{{\rm (odd)}h}\}$ consists of $n/(2h,n)$ elements. The
quotient of these two cardinalities is $(2h,t)/t$. \qed\\

\noindent For future use we make the definition
$r_{g,t}(p):=t_h\rho_{{\rm sgn}(g),(d(g_0)/p),t}((\mathbb Z/p\mathbb Z)^*)$.
Note that $r_{g,t}(p)\in \{0,1,2\}$.
The evaluation of $\rho_{\pm 1,\pm 1,t}(G)$ 
yields the following more precise result for
$r_{g,t}(p)$. (Recall that
$\epsilon_1$ and $\epsilon_2$ are defined in Theorem 1,
respectively Lemma \ref{charactersum}.)
\begin{Lem}
If $g>0$ and $p\equiv 1({\rm mod~}t)$, then
$r_{g,t}(p)=1+\epsilon_2(d(g_0)/p)$. If  $g<0$ and
$p\equiv 1({\rm mod~}2^{1-\epsilon_2}t)$, then
$$r_{g,t}(p)=1+|\epsilon_1|(-1)^{p-1\over 2^{e+1}}\left({d(g_0)\over p}\right).$$
In all other cases $r_{g,t}(p)=0$.
\end{Lem}
Thus, 
\begin{equation}
\label{rhonaarr}
\rho_{{\rm sgn}(g),(d(g_0)/p),t}((\mathbb Z/p\mathbb Z)^*)=r_{g,t}(p)/t_h.
\end{equation}
\indent Using Lemma \ref{kwadstat} and
Lemma \ref{linquad} one easily infers that
$H_{g,t}(x)=M_{g,t}(x)$. Thus, irrespective of the sign
of $g$, we have
\begin{equation}
\label{neededinproof}
H_{g,t}(x)=\sum_{p\le x}\rho_{{\rm sgn}(g),({d(g_0)\over p}),t}
((\mathbb Z/p\mathbb Z)^*)=
{1\over t_h}\sum_{p\le x\atop
p\equiv 1({\rm mod~}t)}r_{g,t}(p)=M_{g,t}(x).
\end{equation}
Using Theorem \ref{surprise1} we see that the
`quadratic' heuristic proposed here is actually asymptotically exact !
The `linear' heuristic, on the other hand, is only asymptotically exact in some cases.\\

\vfil\eject
\section{Equal residual indices}
By inclusion and exclusion it follows that
\begin{equation}
\label{iniexi}
N_{g,t}(x)=\sum_{k=1}^{\infty}\mu(k)R_{g,kt}(x)
\end{equation}
Assuming the error terms to cancel sufficiently, we 
expect from Lemma \ref{algnotheory} that
$$N_{g,t}(x)=A(g,t){x\over \log x}+o\left({x\over \log x}\right).$$
Unfortunately it seems out of reach of present day methods to
prove the cancellation in the error terms. On assuming GRH, however,
the individual error terms involved are all sufficiently small
resulting in a total error term of $o(x/\log x)$, cf.
Theorem \ref{sammie}.

\subsection{Heuristics for equal residual indices}
Just as we used the principle of inclusion and exclusion to study
$N_{g,t}(x)$ in the previous section, we can use it to
set up heuristics for equal residual indices. The analog 
$\sigma_{1,*,t}(G)$ of
the `linear' heuristic $\rho_{1,*,t}(G)$ is given
and evaluated in the next lemma. Note that $\sigma_{1,*,t}(G)$
is the density of elements in $G^h$ having residual index $t$.
\begin{Lem}
\label{inclucond}
Let $G$ be a finite
cyclic group of order $n$ and $t|n$.
We have
$$\sigma_{1,*,t}(G):=
\sum_{d|n/t}\mu(d){|G^h \cap G^{dt}|\over |G^h|}
=\cases{(h,t)\varphi(n/t)/n &if $(n/t,h_t)=1$;\cr
0 &otherwise.}
$$
\end{Lem}
{\it Proof}. By Lemma \ref{intersectie} the sum
under consideration equals
\begin{equation}
\label{zoveelste}
{(h,t)\over t}\sum_{d|n/t}{\mu(d)\over d}
{(h,dt)\over (h,t)}.
\end{equation}
The argument of the latter sum is multiplicative and we
find that it equals zero iff there is a prime divisor
$q$ of $n/t$ satisfying $(h,qt)=q(h,t)$. This
is the case iff $(n/t,h_t)>1$. If $(n/t,h_t)=1$,
then we find that the sum under consideration equals
$t_h^{-1}\prod_{q|n/t}(1-1/q)=(h,t)\varphi(n/t)/n$.
\qed\\

\noindent For a cyclic group $G$ of order $n$
divisible by $t$ let us define
$$\sigma_{\pm 1,\pm 1,t}(G)=\sum_{d|n/t}\mu(d)\rho_{\pm 1,\pm 1,dt}(G).$$
This reduces
to $\sum_{d|(p-1)/t}\mu(d)r_{g,dt}(p)(h,dt)/(dt)$
in case $G=(\mathbb Z/p\mathbb Z)^*$.
The following result holds true (for notational convenience 
we denote $(h,t)\varphi((p-1)/t)/(p-1)$ by $\mu_{g,t}(p)$).
\begin{Lem}
\label{elf}
We have
$\sigma_{{\rm sgn}(g),(d(g_0)/p),t}((\mathbb Z/p\mathbb Z)^*)=w_{g,t}(p)
\mu_{g,t}(p).$
\end{Lem}
{\it Proof}. There are several cases to be considered and we
deal only with a more challenging one: $g<0$ and $2|h_t$
(note that $2|h_t$ is equivalent with $\tau<e$).
If ord$_2(p-1)=\tau$, then
$\sigma_{(d(g_0)/p),t}((\mathbb Z/p\mathbb Z)^*)=0$,
by Lemma \ref{kwadstat}. If ord$_2(p-1)\ge \tau+1$,
that is $p\equiv 1({\rm mod~}2t)$, then by Lemma
\ref{kwadstat},
$\sigma_{(d(g_0)/p),t}((\mathbb Z/p\mathbb Z)^*)$ equals
the sum in (\ref{zoveelste}) but with the divisors $d$ restricted
by ord$_2(p-1)\ge \tau+{\rm ord}_2(d)+1$. This is nothing
but the sum in (\ref{zoveelste}) with $(p-1)/t$ replaced
by $(p-1)/2t$. Thus if $((p-1)/2t,h_t)>1$, then this sum
is zero. If $((p-1)/2t,h_t)=1$, then since $h_t$ is even,
$(p-1)/2t$ is odd and $\varphi((p-1)/2t)=\varphi((p-1)/t)$.
Using this we see that the sum equals
$2(h,t)\varphi((p-1)/2t)/(p-1)=2\mu_{g,t}(p)$.
It follows that if $g<0$ and $2|h_t$, then
$\sigma_{(d(g_0)/p),t}((\mathbb Z/p\mathbb Z)^*)/\mu_{g,t}(p)$
equals 0 if $p\not\equiv 1({\rm mod~}2t)$ and 2
otherwise. These values match with $w_{g,t}(p)$.\\
\indent In the remaining cases sums of the form
(\ref{zoveelste}) appear, but with $d$ restricted to be
even or odd. These sums are easily evaluated. \qed

\begin{Cor}
\label{gevolg}
We have
$\sum_{d|{p-1\over t}}\mu(d)r_{g,dt}(p){(h,dt)\over dt}=w_{g,t}(p)
(h,t){\varphi((p-1)/t)\over p-1}.$
\end{Cor}
The latter corollary expresses $w_{g,t}(p)$ in terms of $r_{g,*}(p)$'s.
It is also possible to express $r_{g,t}(p)$ in terms of $w_{g,*}(p)$'s.
To that end one has to realize that since $\rho_{1,*,t}(G)$ 
and $\sigma_{1,*,t}(G)$
are the fraction
of elements in $G^h$ having residual index divisible by $t$, respectively
equal to $t$, we have $\rho_{1,*,t}(G)=\sum_{d|n/t}\sigma_{1,*,dt}(G)$.
Similarly we have $\rho_{\pm 1,\pm 1,t}(G)=\sum_{d|n/t}
\sigma_{\pm 1,\pm 1,dt}(G)$ and this leads,
on invoking (\ref{rhonaarr}) and Lemma 
\ref{elf}, to the following result.
\begin{Lem}
We have
$r_{g,t}(p)={1\over t_h}\sum_{d|(p-1)/t}w_{g,dt}(p)(h,dt)
{\varphi({p-1\over dt})\over p-1}.$
\end{Lem}
The latter result can be proved also by something akin to
the M\"obius inversion formula:
\begin{Lem}
Let $t$ and $n$ be arbitrary integers with $t|n$ and $\sigma_1$ and $\sigma_2$ be two arithmetic functions,
then $\sum_{d|n/t}\sigma_1(dt)=\sigma_2(t)$ implies
$\sigma_2(t)=\sum_{d|n/t}\mu(d)\sigma_1(dt)$ and vice versa.
\end{Lem}
{\it Proof}. This result
is a particular case of one of Rota's M\"obius inversion formulae
(\cite[Corollary 1, p. 345]{RO}).
If $P$ is a locally finite partially ordered set (whose order
relation is denoted by $\le$) and $r(x)$ is a function on $P$
and $s(x)=\sum_{x\le y\le z}r(y)$, then
$r(x)=\sum_{x\le y\le z}\mu(x,y)s(y)$, where
$\mu(x,y)$ is defined inductively as follows: $\mu(x,x)=1$ for all
$x\in P$. Suppose now that $\mu(x,z)$ has been defined for all $z$ in
the open segment $[x,y)$. Then set
$\mu(x,y)=-\sum_{x\le z<y}\mu(x,z)$.
We apply this with $P$ the partially ordered set of multiples of $t$
dividing $n$, with $x=t$ and $z=n$. On noting that
 $\mu(d,dt)=\mu(d)$, the result follows. \qed\\
\indent Using Lemma \ref{elf} we see that Theorem \ref{main} can
be interpreted as stating that
the `quadratic' heuristic for $N_{g,t}(x)$ is exact up to
order $O(x\log \log x \log^{-2}x)$, under GRH. Indeed, if $N_{g,t}(x)$ tends
to infinity with $x$, then under GRH we have that the `quadratic
heuristic' for $N_{g,t}(x)$ is asymptotically exact.

\vfil\eject
\section{Proof of Theorem 1}
\label{newproof}
In this section we present a proof of Theorem 1 that is rather different
from the one given in \cite{MO}.\\

\noindent {\it Proof of Theorem} 1. Let $C>1$ be arbitrary. The implied constants
below
may depend on $C$, but on $C$ only.
Put $I_1=\sum_{ktd(g_0)\le \log^Cx}\mu(k)M_{g,kt}(x)$
and $I_2=\sum_{ktd(g_0)>\log^C x}\mu(k)M_{g,kt}(x)$. We evaluate
the (finite) sum $I:=I_1+I_2$ in two ways, yielding the proof
on invoking Theorem \ref{sammie}.\\
\indent By Lemma \ref{mtxasymptotic} we have
$$I_1={\rm Li}(x)\sum_{ktd(g_0)\le \log^C x}
{\mu(k)\over [\mathbb Q(\zeta_{kt},g^{1/kt}):\mathbb Q]}
+O\left({x\over \log^C x}\right).$$
Since
$r_{g,t}(p)\le 2$, it
follows by (\ref{neededinproof}) that
$M_{g,t}(x)\le 2h\pi(x;t,1)/t$ and thus $M_{g,t}(x)=0$ for $x>t-1$. 
From this, the latter estimate and the theorem of
Brun-Titchmarsh,
which states that 
the
estimate $\pi(x;t,1)=O(x/(\varphi(t)\log (x/t)))$
holds true
uniformly for $1\le t<x$, we find $I_2=O(hd(g_0)x\log^{-C}x)$.
Using Lemma \ref{wagstaffdegree} we find that
$$\sum_{ktd(g_0)>\log^C x}{\mu(k)\over
[\mathbb Q(\zeta_{kt},g^{1/kt}):\mathbb Q]}
=O\left({hd(g_0)\over \log^{C}x}\right).$$
Combining the latter estimate with those
for $I_1$ and $I_2$ gives
\begin{equation}
\label{eerste}
I=A(g,t){\rm Li}(x)+O(hd(g_0)x\log^{-C}x).
\end{equation}
On the other hand we have, on using
(\ref{neededinproof}) and Corollary \ref{gevolg},
\begin{eqnarray}
I &=&\sum_{k=1}^{\infty}\mu(k)M_{g,kt}(x) =
\sum_{k=1}^{\infty}\mu(k){(h,kt)\over kt}\sum_{p\le x\atop
p\equiv 1({\rm mod~}kt)}r_{g,kt}(p)\nonumber\\
& = &\sum_{p\le x\atop
p\equiv 1({\rm mod~}t)}\sum_{k|{p-1\over t}}\mu(k){(h,kt)\over kt}r_{g,kt}(p)
=(h,t)\sum_{p\le x\atop
p\equiv 1({\rm mod~}t)}w_{g,t}(p){\varphi((p-1)/t)\over p-1}.\nonumber
\end{eqnarray}
Theorem 1 now follows from the
latter equality, (\ref{eerste}) and
Theorem \ref{sammie}. \qed

\medskip\noindent {\footnotesize Korteweg-de Vries Instituut\\
Universiteit van Amsterdam\\
Plantage Muidergracht 24\\
 1018 TV Amsterdam\\
  The Netherlands.\\
e-mail: {\tt moree@science.uva.nl}\\
homepage: 
http://staff.science.uva.nl/$~\tilde{}$moree/}
\end{document}